\documentclass[12pt]{article}
\title{On the ``finitary'' Ramsey's theorem }
\author{Florian Pelupessy\footnote{This research was supported by a JSPS postdoctoral fellowship for foreign researchers}}
\usepackage[english]{babel}
\usepackage{amsmath, amssymb}
\usepackage{url}
\usepackage{xcolor}
\usepackage{hyperref}
\usepackage{tikz}
\usepackage{multicol}
\usepackage{parskip}
\usepackage{makeidx}
  \makeindex

 \usepackage[initials]{amsrefs}

\BibSpec{article}{%
  +{}{\PrintAuthors}				{author} %
  +{,}{ \textit}					{title}
  +{,}{ }							{journal}
  +{}{ \textbf}					{volume}
  +{,}{ }							{note}
  +{}{ \parenthesize}				{date}
  +{:}{ }							{pages}
  +{.}{}							{transition}
  +{}{ }							{review}
  +{\\}{ \url }							{url}
}

\BibSpec{thesis}{%
  +{}{\PrintAuthors}				{author} %
  +{,}{ \textit}					{title}
  +{,}{ }							{journal}
  +{}{ \textbf}					{volume}
  +{,}{ }							{note}
  +{,}{ }							{date}
  +{:}{ }							{pages}
  +{.}{}							{transition}
  +{}{ }							{review}
  +{\\}{ \url }							{url}
}

\usepackage[T1]{fontenc}
\usepackage{libertine}

\newtheorem{theorem}{Theorem}
\newtheorem{lemma}[theorem]{Lemma}
\newtheorem{definition}[theorem]{Definition}

\newtheorem{question}[theorem]{Question}

\newtheorem{corollary}[theorem]{Corollary}

\newcommand{\qed}{
\begin{flushright}
$\Box$
\end{flushright}
}

\newcommand{\RCA}{\mathrm{RCA}_0}
\newcommand{\RCAst}{\mathrm{RCA}_0^{\displaystyle{*}}}
\newcommand{\id}{\mathrm{id}}
\newcommand{\N}{\mathbb{N}}
\newcommand{\WKL}{\mathrm{WKL}_0}
\newcommand{\PA}{\mathrm{PA}}
\newcommand{\ACA}{\mathrm{ACA}_0}
\newcommand{\FRT}{\mathrm{FRT}}
\newcommand{\AS}{\mathrm{AS}}
\newcommand{\RT}{\mathrm{RT}}
\newcommand{\PH}{\mathrm{PH}}
\newcommand{\MD}{\mathrm{MD}}
\newcommand{\WO}{\mathrm{WO}}
\newcommand{\UI}{\mathrm{UI}}
\newcommand{\CF}{\mathrm{CF}}
\newcommand{\AR}{\mathrm{AR}}

\begin{document}
\maketitle
\abstract{
We examine a version of Ramsey's theorem based on Tao, Gaspar and Kohlenbach's ``finitary'' infinite pigeonhole principle. We will show that the ``finitary'' infinite Ramsey's theorem naturally gives rise to statements at the level of the infinite Ramsey's theorem, Friedman's infinite adjacent Ramsey theorem (well-foundedness of certain ordinals up to $\varepsilon_0$), $1$-consistency of theories up to $\PA$ and the finite Ramsey's theorem.
}

\section{Introduction}
This research is inspired by Andreas Weiermann's phase transition programme. The theme of that programme is the following curious phenomenon in first order logic: 

Given a statement $\varphi$ independent of some theory $T$, we can insert a parameter $f\colon \N \rightarrow \N$ in the statement to obtain $\varphi_f$ which may be provable in, or independent of $T$, depending on the parameter value. When one classifies the parameter values $f$ according to the provability of $\varphi_f$ it turns out that, at a threshold value, small changes to $f$ turns $\varphi_f$ from provable in (a weak subtheory of) $T$ to independent of $T$. 

More information on this programme can be found at \cite{weiermannweb}. Our goal in this note is to explore the following question: \emph{What about phase transitions for second order logic?} 

A lazy answer to this question is provided by conservation results, for example: $\ACA$ is conservative over $\PA$, so any phase transition result for $\PA$ is also valid for $\ACA$. However, we may search for more interesting cases in reverse mathematics. Reverse mathematics is the programme, started by Harvey Friedman and, among others, developed by Stephen Simpson, which aims to classify mathematics theorems according to the axioms which are required to prove them. For an introduction to reverse mathematics see \cite{sosoa}. In reverse mathematics we examine \emph{equivalences}. 

Again we may answer our question lazily by restating existing phase transition results, due to the fact that the independent statements used for phase transitions are known to be equivalent to the $1$-consistency of the theory $T$ under consideration. Somewhat less easily, we can also convert existing proofs of these equivalences to show the following: take $\psi_G \equiv \forall f \in G \varphi_f$ and $\alpha$ equal to the proof theoretic ordinal of $T$.
\begin{enumerate}
  \item If $G=\{ f\colon \N \rightarrow \N \}$ then $\psi (G)$ is equivalent to the well-foundedness of $\alpha$.
  \item If $G=\{ f : f \leq \id \}$ then, as stated earlier, $\psi (G)$ is equivalent to the $1$-consistency of the theory $T$.
  \item If $G=\{$constant functions$\}$ then $\psi (G)$ is provable in $\RCA$.
\end{enumerate}
In this note we will examine a more interesting case, where $\psi_G$ has parameter values for which $\psi_G$ is independent of the well-foundedness of $\beta$ for all primitive recursive ordinals $\beta$. 

The starting point is Tao's ``finitary'' pigeonhole principle \cite{taoblog}, which has been extensively studied in \cite{gasparkohlenbach} from the viewpoint of reverse mathematics. We will examine a ``finitary'' version of Ramsey's theorem which is a generalisation of Tao's pigeonhole principle. 
\begin{definition}[$\AS$]
A function $F\colon \{$(codes of) finite subsets of $\N\} \rightarrow \N$ is asymptotically stable if for every sequence $X_0 \subseteq X_1 \subseteq X_2 \dots$ of finite sets, there exists $i$ such that $F(X_j)=F(X_i)$ for all $j \geq i$.
\end{definition}
This definition of $\AS$ is modified from \cite{taoblog}. Roughly speaking, $|X| \geq F(X)$ can be interpreted as `the finite set $X$ is large'. $\AS$ would then be the set of possible manners in which to define `large'. 
\begin{definition}[$\FRT^k_d$] ~{} \\
For every $F \in \AS$ there exists $R$ such that for all  $C\colon [0,R]^d \rightarrow k$ there exists $C$-homogeneous $H$ of size $> F(H)$.
\end{definition}
\begin{definition}
$\FRT_d$ is the statement $\forall k.\FRT^k_d$. $\FRT$ is the statement $\forall d,k . \FRT^k_d$.
\end{definition}
\begin{definition}[$\RT^k_d$]  ~{} \\
For every $C\colon [\N]^d \rightarrow k$ there exists an infinite $C$-homogeneous set.
\end{definition}
One can view $\FRT$ as the collection of all finite versions of $\RT$, similar to the familiar finite Ramsey's theorem. We will show that, as is shown for the case $d=1$ in\cite{gasparkohlenbach}, $\FRT^k_d$ is quivalent to $\RT^k_d$ over $\WKL$. \\

\noindent Notice the following: \\

\noindent If, in $\FRT$, we replace $\AS$ with the set of constant functions: 

\begin{definition}[$\CF$]
\[
\exists m. F=m,
\]
\end{definition}
the resulting theorem becomes simply the finite Ramsey's theorem. \\  

\noindent If we replace $\AS$ with the following:
\begin{definition}[$\UI$]
\[
\exists m \forall X . F(X)\leq \max \{ \min X , m \},
\]
\end{definition}
then the resulting theorem is the Paris--Harrington principle, which, for dimension $d+1$ is equivalent to the $1$-consistency of $\mathrm{I}\Sigma_d$. It is equivalent to $1$-consistency of $\PA$ for unrestricted dimensions. 
\begin{definition}
$\FRT^k_d(G)$ is the statement obtained from $\FRT^k_d$ by replacing $F \in \AS$ with $F\in G$. $\FRT_d(G)$ is the statement $\forall k.\FRT^k_d(G)$. $\FRT(G)$ is the statement $\forall d,k . \FRT^k_d(G)$.
\end{definition}
One obvious question is whether there are properties $G$ such that the strength of $\FRT(G)$ lies strictly between $\FRT(\UI )$ and $\FRT(\AS )$. We will show that this is the case for:
\begin{definition}[$\MD$]
\[
\forall X,Y . \min X=\min Y \rightarrow F(X)=F(Y).
\]
\end{definition}
Because this latter version has connections with Friedman's adjacent Ramsey theorem we conclude with determining the level-by-level strength of the adjacent Ramsey theorem.

\section{$\FRT$}
We assume familiarity with reverse mathematics, primitive recursion, $\RCA$, $\WKL$ and Ramsey's theorem as in Chapters II and IV in  \cite{sosoa}. Please note that for finite set $X$ we also use $X$ to denote its code. \\

\noindent The main theorem in this section is:

\begin{theorem} \label{thm:mainfrt} ~{}
\begin{enumerate} 
  \item[(a)] $\RCA \vdash \FRT^k_d \rightarrow \RT^k_d$,
  \item[(b)] $\WKL \vdash \RT^k_d \rightarrow \FRT^k_d$.
\end{enumerate}
\end{theorem}
We will make use of:
\begin{lemma}\label{lemma:primrec} The following are primitive recursive:
\begin{enumerate}
  \item the relation $\{ (x, X) : x \in X\}$,
  \item the relation $\{ (X, Y) : X \subseteq Y \}$,
  \item the relation $\{ (X, C) : X$ is $C$-homogeneous$\}$ and
  \item the function $(x, C) \mapsto C(x)$ for finite functions $C$.
\end{enumerate}
\end{lemma}
\emph{Proof:} Exercise for the reader.
\qed
\emph{Proof of Theorem~\ref{thm:mainfrt} (a):} We adapt the proof of the case $d=1$ from \cite{gasparkohlenbach}. Please notice the extra steps needed to deal with the modified definition of $\AS$.

In $\RCA$, we show $\neg \RT^k_d \rightarrow \neg \FRT^k_d$. Suppose $C\colon [\N]^d \rightarrow k$ is a colouring such that every $C$-homogeneous set has finite size.  Define the following $F$ primitive recursively:
\[
F(X) = 
\left\{ 
\begin{array}{ll}
|X|+1 & \textrm{ if $X$ is $C$-homogeneous}, \\
0 & \textrm{otherwise}.
\end{array}
\right.
\]
\emph{Claim 1:} $F \in \AS$. Take $X_0 \subseteq X_1 \subseteq \dots$. Examine the $\Sigma^0_1$ formula:
\[
\varphi(n) \equiv \exists i (n \in X_i).
\]
By Lemma II.3.7 of \cite{sosoa} $\{n : \varphi(n)\}$ is finite or there exists a one-to-one function $f$ such that
\[
\forall n [\varphi(n) \leftrightarrow \exists m (f(m)=n)].
\]

If $\{n : \varphi(n)\}$ is finite then there exists $i$ with $F(X_i)=F(X)$ and we are finished with the claim, so assume the latter case.

We will show that there exists an infinite set $X$ such that $n \in X \rightarrow \varphi(n)$ (hence $X$ is a subset of the possibly nonexistent $\bigcup X_i$). This is sufficient, because then $\forall i \exists j>i F(X_j) \neq F(X_i)$ implies $X$ is $C$-homogeneous. We show this by translating a rather common exercise from computability theory to our context: Given an infinite recursively enumerable set, show that it contains an infinite decidable subset.
\\

\noindent Take $\Sigma^0_1$ formula:
\[
\phi(n) \equiv \exists m [f(m) \geq n \wedge f(\mu x \leq m. f(x) \geq n)=n].
\]
and $\Pi^0_1$ formula:
\[
\psi(n) \equiv \forall m [f(m) \geq n \rightarrow f(\mu x \leq m. f(x) \geq n)= n].
\]

These two formulas are equivalent by unboundedness of $f$, so by $\Delta^0_1$-comprehension the infinite set $X=\{ n : \phi(n)\}$ exists. This finishes the proof of claim 1. \\

\noindent \emph{Claim 2:} $F$ is a counterexample for $\FRT^k_d$. Take arbitrary $R$, Define $D=C$ restricted to $[0,R]^d$. By definition of $F$ any $D$-homogeneous set $H$ has size $<F(H)$, ending the proof of claim 2 and part (a) of the theorem. \\

\noindent \emph{Proof of Theorem~\ref{thm:mainfrt} (b):} We use a compactness proof which involves K\"onig's lemma. However, we take care that the application of K\"onig's lemma uses only the bounded version (hence we reason in $\WKL$ by Lemma IV.1.4 in \cite{sosoa}).

Assume $\neg \FRT^k_d$, hence there exists $F \in \AS$ such that for all $R$ there exists $C\colon [0,R]^d \rightarrow k$ for which every $C$-homogeneous set $H\subseteq [0,R]$ has size $\leq F(H)$. Enumerate such colourings with $\{C_{R,i}\}_{i\leq n_R}$. Notice that the codes of these colourings can be bounded by some function which is primitive recursive in $d,k,R$. We define the following bounded (by previous remark) and infinite tree: 
\[
T = \{ \langle C_{1, i_1} , \dots , C_{R, i_R} \rangle : C_{1,i_1} \subseteq \dots \subseteq C_{R, i_R} \}.
\]

Take the colourings $D_1 \subseteq D_2 \subseteq \dots$ from the infinite path in $T$, which exists due to bounded K\"onig's lemma. Define $D\colon [\N]^d \rightarrow k$ as follows:
\[
D(x) = D_{\max x} (x).
\]
\emph{Claim:} $D$ is a counterexample for $\RT^k_d$. Assume $H$ is $D$-homogeneous. By construction of $T$ and $D=\bigcup D_i$, the size of $H_i =H \cap [0, i]$ is less than or equal to $F(H_i)$ for every $i$. Note that $H_1 \subseteq H_2 \subseteq H_3 \subseteq \dots$ and $H=\bigcup H_i$, so (by $F \in \AS$) there exists $i$ such that $F(H_j)=F(H_i)$ for all $j \geq i$, hence $H$ is finite. This ends the proof of the claim and part (b) of the theorem. 
\qed
\begin{question}
Is $\WKL$ required in part (b) of this theorem? Notice that $\WKL$ is not required  for $d \geq 3$.
\end{question}

\section{Restriction to the minimally dependent}

We assume basic familiarity with ordinals up to $\varepsilon_0$ and their cantor normal forms. 
\begin{definition}
$\omega_0=1$ and $\omega_{n+1} = \omega^{\omega_n}$.
\end{definition}
\begin{definition}[$\WO(\alpha)$]
Every infinite sequence $\alpha_0, \alpha_1, \dots$ below $\alpha$ has $i<j$ such that $\alpha_i \leq \alpha_j$.
\end{definition}
The main theorem in this section is:
\begin{theorem} \label{thm:parph}
$\RCA \vdash \WO (\omega_d) \leftrightarrow \FRT_d (\MD )$
\end{theorem}
Observe first that $\FRT_d(\MD )$ is equivalent to $\forall f\colon \N \rightarrow \N. \PH^d_f$.
\begin{definition}[$\PH^d_{f}$]
For all $a$ there exists $R$ such that for all $C\colon [a,R]^d \rightarrow k$ there exists $C$-homogeneous $H$ of size $f(\min H)$.
\end{definition}
\emph{Proof of Theorem~\ref{thm:parph}:} `$\rightarrow$' in Subsection~\ref{subsection:upper} \\
`$\leftarrow$'  in Subsection~\ref{subsection:lower}. 
\qed
\subsection{Lower bound} \label{subsection:lower}
We modify the proof of $\PH_\id \rightarrow$ $\mathrm{Tot}(H_{\varepsilon_0})$ from \cite{pelupessyfriedman}. The proof below consist mostly of recalling the necessary definitions and lemmas, where the final step is modified to fit our new situation. We skip the proofs when they are unchanged from the original.
\begin{definition} 
Given $\alpha=\omega^{\alpha_1} \cdot a_1 + \dots + \omega^{\alpha_n} \cdot a_n$ and $\beta=\omega^{\beta_1} \cdot b_1+ \dots + \omega^{\beta_m} \cdot b_m$, with the $a_i,b_i$ positive integers, $\alpha_1 > \dots > \alpha_n$ and $\beta_1 > \dots > \beta_m$ we define: 
\begin{enumerate}
\item The comparison position $\mathrm{CP}(\alpha , \beta)$ is the smallest $i$ such that $\omega^{\alpha_i} \cdot a_i \neq \omega^{\beta_i} \cdot b_i$ if such an $i$ exists, zero otherwise.
\item The comparison coefficient $\mathrm{CC}(\alpha , \beta)$ is $a_{\mathrm{CP}(\alpha, \beta)}$, where $a_0=0$.
\item The comparison exponent $\mathrm{CE}(\alpha , \beta )$ is $\alpha_{\mathrm{CP}(\alpha, \beta)}$, where $\alpha_0=0$.
\end{enumerate}
Define the maximal position $\mathrm{MP}$ and maximal coefficient $\mathrm{MC}$ by induction on $\alpha$ as follows: $\mathrm{MP}(0)=1$ and $\mathrm{MC}(0)=0$. Given $\alpha=\omega^{\alpha_1} \cdot a_1 + \dots + \omega^{\alpha_n} \cdot a_n>0$, with the $a_i$ positive integers and $\alpha_1 > \dots > \alpha_n$, define:
\begin{enumerate}
\item[(4)] $\mathrm{MP}(\alpha )=\max \{ n, \mathrm{MP}(\alpha_i) \}$.
\item[(5)] $\mathrm{MC}(\alpha)=\max \{ a_i , \mathrm{MC}(\alpha_i) \}$.
\end{enumerate}  
\end{definition}
\begin{lemma}
We have:
\begin{enumerate} \label{lessimpliesless}
\item $\mathrm{CP}(\alpha , \beta) \leq \mathrm{MP}(\alpha )$.
\item $\mathrm{CC}(\alpha , \beta) \leq \mathrm{MC}(\alpha )$.
\item $\mathrm{MP}(\alpha_i) \leq \mathrm{MP}(\alpha)$ and $\mathrm{MC}(\alpha_i) \leq \mathrm{MC}(\alpha)$.
\item $\mathrm{CP}(\alpha , \beta) \leq \mathrm{CP}(\beta, \gamma ) \wedge \mathrm{CE}(\alpha , \beta) \leq \mathrm{CE}(\beta, \gamma ) \wedge \mathrm{CC}(\alpha , \beta) \leq \mathrm{CC}(\beta, \gamma ) \Rightarrow \alpha \leq \beta$.
\end{enumerate}
\end{lemma}
\begin{definition} \index{$F_d^l$} \label{def:function}
Let $l,d,n$ be nonnegative integers. Define $\omega_0(l)=l$ and $\omega_{n+1}(l)=\omega^{\omega_n(l)}$. Define $F_d^l\colon \omega_{d}(l+1)^d \rightarrow \mathbb{N}^{2d+l-1}$ by recursion on $d$:
\begin{enumerate}
\item Given $\alpha=\omega^l\cdot n_l + \dots + \omega^0 \cdot n_0$, define $F_1^l(\alpha)=(n_l, \dots , n_0)$.
\item $F_{d+1}^l (\alpha_1 ,\dots ,\alpha_{d+1})=$ \\
\ \ $( \mathrm{CP}(\alpha_1, \alpha_2), \mathrm{CC}(\alpha_1, \alpha_2), F_d^l (\mathrm{CE} (\alpha_1, \alpha_2), \dots , \mathrm{CE}(\alpha_d, \alpha_{d+1})))$.
\end{enumerate}
\end{definition}
\begin{lemma} \label{lemma:tussen2}
$F_d^l(\alpha_1, \dots, \alpha_d) \leq F_d^l(\alpha_2, \dots , \alpha_{d+1}) \Rightarrow \alpha_1 \leq \alpha_2$.
\end{lemma}
\begin{lemma} \label{lemma:tussen}
$F_d^l(\alpha_1, \dots, \alpha_d) \leq \max \{ \mathrm{MC}(\alpha_1), \mathrm{MP}(\alpha_1) \}$.
\end{lemma} 
We are finally ready to finish the proof the lower bound of Theorem~\ref{thm:parph}. The following lemma is where the proof from \cite{pelupessyfriedman} is modified:
\begin{lemma}
$\RCA \vdash \forall f . \PH_f^d \rightarrow \mathrm{WO}(\omega_d)$, where $\PH_f^d$ is $\PH_f$ with fixed dimension $d$.
\end{lemma}
Given infinite sequence $\alpha_0, \alpha_1, \alpha_2, \dots$ below $\omega_{d}(l+1)$ take 
\[
f(i)=\max\{\mathrm{CC}(\alpha_i), \mathrm{CP}(\alpha_i)\}+d+2
\]
and $R$ from $\PH_f$ in dimension $d+1$, $a=0$ and $c=2d+l$. Define colouring $C\colon [0, R]^{d+1}\rightarrow  [0,2d+l]$:
\[
C(x_1, \dots , x_{d+1})=\left\{ 
\begin{array}{ll}
0 & \textrm{if $F_d^l(\alpha_{x_1}, \dots, \alpha_{x_d }) \leq$} \\
   & \textrm{\ \ \   $F_d^l(\alpha_{x_2 }, \dots, \alpha_{x_{d+1} })$}, \\
i  & \textrm{otherwise},
\end{array}
\right.
\]
where $i$ is the least such that:
\[
(F_d^l(\alpha_{x_1 }, \dots, \alpha_{x_d }))_i > (F_d^l(\alpha_{x_2}, \dots, \alpha_{x_{d+1} }))_i.
\]
Observe that $(F_d^l(\alpha_{x_1 }, \dots, \alpha_{x_d }))_i \leq  \max\{\mathrm{CC}(\alpha_{x_1}), \mathrm{CP}(\alpha_{x_1})\}$ (this is a consequence  of Lemma~\ref{lemma:tussen} ). Take homogeneous $H$ of size $f(\min H)$from $\PH_f$. If the value of $C$ on $[H]^{d+1}$ is $i>0$ we can obtain a decending sequence of natural numbers below $f(\min H) - d - 2$ of length $f(\min H) -d$, which is impossible. Hence the value of $C$ is $0$, delivering $\alpha_{x_1} \leq \alpha_{x_2}$. 
\qed
\subsection{Upper bound} \label{subsection:upper}
We use the upper bounds result from Section 6 in \cite{ketonensolovay}, observing that, mostly thanks to the formalisation of large parts in $\mathrm{I}\Sigma_1$  in Section II.3 in \cite{hajekpudlak}, the proofs are within $\RCA+ \mathrm{WO}(\omega_d)$. Alternatively, one can use Corollary~15 from \cite{pelupessyKS}, which states that the theorem in question is provable in $\RCA$.

A similar version, called relativised Paris--Harrington for $d=2$ has also been studied by Kreuzer and Yokoyama in \cite{kreuzeryokoyama}.
\begin{definition}
$A=\{a_0 < \dots < a_b \}$ is $\alpha$-large if $\alpha[a_0] \dots [a_b]=0$, where $\alpha [.]$ denotes the canonical fundamental sequences for ordinals below $\varepsilon_0$.
\end{definition}
\begin{lemma}
$\RCA$ proves the following:  if $\WO (\omega_d)$ then  for every strictly increasing  $f\colon \N \rightarrow \N$, $a \in \N$, $\alpha < \omega_d$ there exists $\alpha$-large set $\{ f(a) , f(a+1) ,\dots , f(b) \}$.
\end{lemma}
\emph{Proof:} Define the following descending sequence of ordinals: $\alpha_0=\alpha$ and:
\[
\alpha_{i+1}=\alpha_i [f(i)].
\]
By well-foundedness of $\omega_d$ this sequence reaches zero, delivering the desired $\alpha$-large set.
\qed
Assume without loss of generality, that $f$ is strictly increasing and $>3$. By $\WO (\omega_d)$ there exists $\omega_{d-1} (c+5)$-large set $A=\{ f(a) , \dots , f(b) \}$. We claim that $R=b$ witnesses $\PH^d_f$: Take colouring $C \colon [a,R]^d \rightarrow c$, define $D\colon [A]^d \rightarrow c$ as follows:
\[
D(x_1, \dots , x_d) = C( f^{-1} (x_1) , \dots , f^{-1} (x_d)).
\]
By Theorem 6.7 from \cite{ketonensolovay} or Corollary~15 from \cite{pelupessyKS} there exists $D$-homogeneous $X$ with size $\min X$. Then $H= \{ f^{-1}(x): x\in X\}$ is $C$-homogeneous and of size $f(\min H)$. This ends the proof of Theorem~\ref{thm:parph}.
\qed

\section{$\FRT$ and adjacent Ramsey}
\begin{definition}
For $r$-tuples $a \leq b$ denotes the coordinatewise ordering: 
\[
a \leq b \Leftrightarrow (a)_1 \leq (b)_1 \wedge \dots \wedge (a)_r \leq (b)_r.
\]
\end{definition}
\begin{definition}[$\AR_d$]
For every $C \colon \N^d \rightarrow \N^r$ there exist $x_1 < \dots < x_{d+1}$ such that $C(x_1, \dots , x_d) \leq C(x_2 , \dots x_{d+1})$.
\end{definition}
\begin{definition}
$\AR$ denotes $\forall d . \AR_d$.
\end{definition}
In this section we will show that:
\begin{theorem}
$\RCA \vdash \WO (\omega_{d+1}) \leftrightarrow \AR_d$
\end{theorem}
\emph{Proof:} `$\leftarrow$': We use $F^l_d$ from \ref{subsection:lower}. Given sequence of ordinals $\omega_{d}(l+1)> \alpha_0, \alpha_1, \dots$ define:
\[
C(x_1, \dots , x_d) = F^l_d (\alpha_{x_1} , \dots , \alpha_{x_d}).
\]
By $\AR_d$ there exist $x_1 < \dots < x_{d+1}$ with $C(x_1, \dots , x_d) \leq C(x_2, \dots , x_{d+1})$, which by Lemma~\ref{lemma:tussen2} deliver $\alpha_{x_1} \leq \alpha_{x_2}$. \\

\noindent `$\rightarrow$': By Theorem~\ref{thm:parph} it is sufficient to show that $\forall f . \PH^{d+1}_f \rightarrow \AR_d$.  For this it is sufficient to simply note that the proof of $\PH_{d+1} \rightarrow \AR_d$ from \cite{pelupessyfriedman} (please note the difference in $\AR$ as defined there) works fine when relative to the function
\[
f(x) = \max_{ y \in \{ 0, \dots , x \}^d } C(y).
\]
Replace the strong adjacent Paris--Harrington principle with a version relative to $f$:
\begin{definition}[$\mathrm{SAPH}^d_f$]
For every $c,k.m$ there exists an $R$ such that for every colouring $C \colon [m, \dots , R]^d \rightarrow [0,c]$ there exists $C$-homogeneous $H=\{ h_1 < h_2 < \dots\}$ of size $f(h_k)$.
\end{definition}
Then $\forall f. \PH^{d+1}_f \rightarrow \forall f. \mathrm{SAPH}^{d+1}_f \rightarrow \AR_d$ by copying the proofs of Theorems 3.4 and 3.5 from  \cite{pelupessyfriedman}. 
\qed

\section{Conclusions}

$\RCA$ proves the following: \\
\begin{tabular}{ l c l c l}
$\FRT$ 				& $\leftrightarrow$ &	$\RT$ \\
$\FRT^k_d$ 			& $\leftarrow$	    & $\RT^k_d$ for (d>2) \\
$\FRT^k_d$			& $\rightarrow$	    & $\RT^k_d$ \\
$\FRT (\MD)$ 		& $\leftrightarrow$ & $\AR$ & $\leftrightarrow$ & $\WO (\varepsilon_0)$ \\
$\FRT_{d+1} (\MD)$	& $\leftrightarrow$ & $ \AR_d$ & $\leftrightarrow$ &$\WO (\omega_{d+1})$ \\
$\FRT (\UI)$			& $\leftrightarrow$ & $1$-consistency of $\PA$ \\
$\FRT_{d+1} (\UI)$	& $\leftrightarrow$ & $1$-consistency of $\mathrm{I}\Sigma_d$ \\
$\FRT (\CF)$ \\
\end{tabular}
\\
The last three of those lines are true because $\FRT_d (\UI)$ is equivalent to $\PH^d_{\id}$, so the equivalence to $1$-consistency is the classic Paris--Harrington result from \cite{parisharrington}. \\

\noindent Furthermore, $\WKL \vdash \RT^k_d \rightarrow \FRT^k_d$. 

\begin{corollary} Over $\RCA$:
\[
\FRT (\CF ) < \FRT (\UI ) < \FRT (\MD ) < \FRT (\AS ). 
\]
\end{corollary}

\begin{question}
Do the same implications hold for $\RCA^*$ and, where $\WKL$ is used, in $\WKL^*$? 
\end{question}

\subsubsection*{Acknowledgements:} 

\noindent The author thanks Kazuyuki Tanaka and Keita Yokoyama for discussions on Reverse Mathematics and Ramsey theorem variants. \\

\noindent Additionally, the author is grateful to Makoto Fujiwara for pointing out the modified definition of $\AS$, which has an existential set quantifier removed compared to the definition in \cite{gasparkohlenbach}.


\begin{thebibliography}{14}


\bib{pelupessyfriedman}{article}{
 author={Friedman, Harvey}
 author={Pelupessy, Florian}
 title={Independence of Ramsey theorem variants using $\varepsilon_0$}
 journal={Proc. Amer. Math. Soc. , in press}
}

\bib{gasparkohlenbach}{article}{
  author={Gaspar, Jaime},
  author={Kohlenbach, Ulrich},
  title={On Taos ``finitary''  infinite pigeonhole principle},
  journal={},
}

\bib{ketonensolovay}{article}{
  author={Ketonen, Jussi}
  author={Solovay, Robert}
  title={Rapidly growing Ramsey functions}
  journal={Annals of Mathematics}
}

\bib{kreuzeryokoyama}{article}{
  author={Kreuzer, Alexander}
  author={Yokoyama, Keita}
  title={On principles between $\Sigma_1$- and $\Sigma_2$-induction and monotone enumerations}
}


\bib{parisharrington}{article}{
  author={Paris, Jeff},
  author={Harrington, Leo},
  title={A Mathematical Incompleteness in Peano Arithmetic},
  journal={in Handbook for Mathematical Logic (Ed. J. Barwise) Amsterdam, Netherlands: North-Holland, 1977}
}

\bib{pelupessyKS}{article}{
  author={Pelupessy, Florian}
  title={On $\alpha$-largeness and the Paris--Harrington principle in $\RCA$ and $\RCAst$}
  journal={\href{https://arxiv.org/abs/1611.08988}{arXiv:1611.08988}}
}

\bib{hajekpudlak}{book}{
  author={H\'ajek, Petr}
  author={Pudl\'ak, Pavel}
  title={Metamathematics of First-order Arithmetic}
  series={Perspectives in Mathematical Logic}
  volume={3}
  year={1998}
  publisher={Berlin: Springer-Verlag}
}

\bib{sosoa}{book}{
  author={Simpson, Stephen},
  title={Subsystems of second order arithmetic},
}

\bib{taoblog}{article}{
  author={Tao, Terrence}
  title={Blog post on Soft analysis, hard analysis}
  url={https://terrytao.wordpress.com/2007/05/23/soft-analysis-hard-analysis-and-the-finite-convergence-principle/}
}

\bib{weiermannweb}{article}{
  author={Weiermann, Andreas},
  title={Webpage on phase transitions},
  url={http://cage.ugent.be/~weierman//phase.html},
}

\end{thebibliography}
\end{document}